\documentclass{elsarticle}
\usepackage{amsmath}
\usepackage{amssymb}
\usepackage{graphics}
\usepackage{epstopdf}
\usepackage{float}
\newtheorem{theorem}{Theorem}[section]
\newtheorem{lemma}{Lemma}[section]
\newtheorem{definition}{Definition}[section]
\usepackage[super]{nth}

\begin{document}

\begin{frontmatter}

\title{Tikhonov Theorem for Differential Equations with Singular Impulses}

\author[a]{Marat ~Akhmet\corref{cor}}
\ead{marat@metu.edu.tr}

\author[a]{Sabahattin ~\c{C}a\u{g}}
\ead{sabahattincag@gmail.com}

\cortext[cor]{Corresponding author}

\address[a]{Department of Mathematics, Middle East Technical University, Ankara, Turkey.}

\begin{abstract}
The   paper   considers impulsive systems with   singularities.  The main novelty of the   present  research  is that impulses  (impulsive functions)    are   singular.  This is beside singularity of differential equations. The Lyapunov second method is applied to proof the main theorems.  Illustrative examples with simulations are given to support the theoretical results.
\end{abstract}


\end{frontmatter}

%\linenumbers

\section{Introduction}
The singularly perturbed differential equations  arise in various fields  of chemical kinetics \cite{segel}, mathematical biology \cite{hek,owen}, fluid dynamics \cite{fluid} and in a variety models for control theory \cite{kokotovic,gondal}. These problems depend on a small positive parameter such that the solution varies rapidly in some regions and varies slowly in other regions. 

The contribution of our work relates to  a new Tikhonov
theorem for  singularly perturbed impulsive systems. This theorem express the limiting behavior of solutions of the singularly perturbed system. It is a powerful instrument for analysis of singular perturbation problems. It has been studied for many types of differential equations; partial differential equations \cite{Kadalbajoo}, singularly perturbed differential inclusions \cite{Veliov}, functional-differential inclusions \cite{Tzanko}, discontinuous differential equations \cite{Chen,Chen2,Chen3,simeonov1,simeonov2}.

Impulse effects exist in a wide diversity of evolutionary processes that exhibit abrupt changes in their states \cite{Akhmet2010, Akhmet2011,akhmet2015}. In many systems, in addition to singular perturbation, there are also impulse effects \cite{Chen,Chen2,Chen3,simeonov1,simeonov2}.
Chen et al. \cite{Chen3} derived a sufficient condition that guarantees robust exponential stability for sufficiently small singular perturbation parameter by applying the Lyapunov function method and using a two-time scale comparison principle. In \cite{simeonov1,simeonov2}, authors  proposed Lyapunov function method to set up the exponential stability criteria for singularly perturbed impulsive systems. This method can be efficiently used to overcome the impulsive perturbation such that the stability of the original system can be ensured. In \cite{Chen},  Lyapunov function method was further extended to study the exponential stability of singularly perturbed stochastic time-delay systems with impulse effect.  The results in \cite{Chen,simeonov1,simeonov2} only guarantee the systems under consideration to be exponentially stable for sufficiently small positive parameter.

Various types of singular perturbation problems are discussed in many books \cite{Tikhonov,bainov,vasileva,Omalley1991,verhulst}. Consider the following model of singularly perturbed differential equation 
%	\begin{equation} \label{tikhonov}
%	\mu \frac{dz}{dt} = F(z),
%	\end{equation} 
%where $\mu$ is a small positive real number. Then using this system, they study the following  well known model of the singularly perturbed problem 
\begin{equation}\label{tikhnov2}
\begin{split}
\mu \dot{z}=f(z,y,t),\\
\dot{y}=g(z,y,t),
\end{split}
\end{equation}
where $\mu$ is a small positive real number. In the literature, the results based on this system are known as  Tikhonov Theorems \cite{Omalley1991,verhulst,tikhonov1952}. 
Bainov and Covachev \cite{bainov} first extended the impulsive analogous of Tikhonov Theorem concerning system \eqref{tikhnov2} in the form of
 \begin{equation}\label{bainov}
 	\begin{split}
 		\mu \dot{z}=f(z,y,t), \quad \Delta z|_{t=t_i}=I_i(y(t_i)),\\
 		\dot{y}=g(z,y,t), \quad \Delta y|_{t=t_i}=J_i(y(t_i)),
 	\end{split}
 \end{equation}
 where $i=1,2,...,p$ and $0<t_1<t_2<...<t_p<T.$  However, only the differential equation in their problem is singularly perturbed.
 
 In this study,  we consider differential equations where impulses are also singularly perturbed which are different than in \cite{bainov}. The following system is our focus of discussion
 \begin{equation} \label{eq1-general}
 \begin{split}
 \mu \frac{dz}{dt} &= F(z,y,t) ,\quad \mu \Delta z|_{t=\theta_{i}}=I(z,y,\mu)\\
 \frac{dy}{dt} &= f(z,y,t) , \quad \Delta y|_{t=\eta_{j}}=J(z,y),
 \end{split}
 \end{equation} 
 where $z,F$ and $I$ are m-dimensional vector valued functions, $y,f$ and $J$ are n-dimensional vector valued functions, $0<\theta_1<\theta_2<\dots<\theta_p<T,$ $\theta_i,i=1,2,\dots,p,$ and $\eta_j,j=1,2,\dots,k,$ are distinct discontinuity moments in $(0,T).$ 
 
 The main novelty of this paper is the extension of Tikhonov Theorem such that system  \eqref{eq1-general} has  the small parameter in impulse function, the discontinuity moments are different for each dependent variables. The singularity in the impulsive part of the system can be treated through perturbation theory methods. 
%%%%%%%%%%%%%%%%%%%%%%%%%%%%%%%%%%%%%%%%%%%%%%%%%%%%%%%%
%%   
%%%%%%%%%%%%%%%%%%%%%%%%%%%%%%%%%%%%%%%%%%%%%%%%%%%%%%%
\subsection{Preliminaries}
Let us describe generally the definition of singularity. Consider
\begin{itemize}
	\item Problem $P(\mu)$: the problem with small parameter $\mu,$
	\item Problem $P(0)$: the reduced (degenerated) problem.
\end{itemize}
The problem $P(0)$ is a simplified model of  $P(\mu)$ taking $\mu=0.$ Denote the solution of $P(0)$ by $z(t,0)$ and the solution of $P(\mu)$ by $z(t,\mu).$
\begin{definition}\cite{vasileva}
	$P(\mu)$ is called regularly perturbed problem in a domain $D$ if 
	$$\sup_D{\|z(t,\mu)-z(t,0)\|}\to 0 \text{ when } \mu \to 0.$$
	Otherwise, it is called singularly perturbed problem.
\end{definition}
It follows from the definition that for a regularly perturbed problem the solution $z(t,0)$ of $P(0)$ will be close to the solution $z(t,\mu)$ of $P(\mu)$ in the entire domain $D$ for all sufficiently small $\mu$. However, if the problem $P(\mu)$ is singularly perturbed, then  $z(t,\mu)$ will not be close to  $z(t,0)$ for all small $\mu$ at least in some part of domain $D.$
\section{A Particular Case of the Main Theorem}
Before carrying out our main investigation, let us consider a particular case of the main theorem. This case is useful by its geometric clarity. We introduce the following problem
  \begin{equation} \label{eq1}
  \begin{split}
  &\mu \frac{dz}{dt} = F(z) , \\
  &\mu \Delta z|_{t=\theta_{i}}=I(z,\mu), 
  \end{split}
  \end{equation} 
  with $z(0,\mu)=z_{0}$, where $z  \in \mathbb{R}^{m}$, $t \in [0, T],$ $F(z)$ is a continuously differentiable function on D  and $I(z,\mu)$ is a continuous function for $(z,\mu)\in D\times[0,1]$ ,  D is the domain $D=\{0\leq t\leq T, \| z \|<d \},$  $\theta_i$s are defined above.
  
The parameter in a the impulsive equation makes it possible that $\frac{I(z,\mu)}{\mu}$ blow up at impulse moments as $\mu \to 0.$ This is why, a deep analysis and convenient conditions for the limiting processes with $\mu \to 0$ have to be researched.
\subsection{Singularity with a Single Layer}
Let us take $\mu=0$ in \eqref{eq1}.  Then, one has $0=F(z)=I(z,0).$  It is degenerate system  since its order is less than the order of \eqref{eq1}. Consider an isolated real root $z=\varphi$ of $F(z)=0$ and $I(z,0)=0.$ %Assume that all the roots $z=\varphi$ of $F(z)=0$ such that $I(\varphi,0)=0$ are real and isolated in $\bar{D}$. It is necessary to choose one of the roots in order to guarantee that the solution $z(t,\mu)$ of \eqref{eq1} with $z_0$ is close to that root. 
%The rule for choosing the root will be stated.

Now,  introduce a new variable $\tau=\frac{t}{\mu}$ and $x=z-\varphi$ for the first equation in \eqref{ex1} to obtain 
  \begin{equation}\label{eq2}
	\frac{dx}{d\tau} = F(x+\varphi). 
  \end{equation}
The following condition is required. 
\begin{itemize}
	\item[(C1)] Suppose that there is a positive definite function $V(x)$ such that $V(0)=0$ and  whose derivative with respect to $\tau$ along system \eqref{eq2} is negative definite.
\end{itemize}
This condition implies that the zero solution of \eqref{eq2} is uniformly asymptotically stable.
Moreover, for the impulsive function we need the following condition.
\begin{itemize}
	\item[(C2)] \[ \lim_{(z,\mu)\to (\varphi,0)}\frac{I(z,\mu)}{\mu}= 0,\]
\end{itemize}
which prevents impulsive function to blow up as the parameter $\mu$ decays to zero. This condition is the counterpart of (C1) considering impulsive function. Condition (C2) plays similar role to condition (C1) in the proof of the next theorem.

\begin{theorem} \label{thmlyapunov1}
 Suppose that conditions (C1) and (C2) are true. If the initial value $z_0$ is located in the domain of attraction of the root $\varphi,$ then  solution $z(t,\mu)$ of \eqref{eq1} with $z(0,\mu)=z_0$ exists on $[0,T]$ and it is satisfies the limit
\begin{equation}\label{lim1}
	\lim_{\mu \to 0}z(t,\mu)=\varphi \quad \text{for} \quad 0<t\leq T.
\end{equation}
\end{theorem}
 
\textbf{Proof.} In this proof, we will follow the idea of the proof of  \cite[Theorem 7.3]{Tikhonov}. Consider the first interval $[0,\theta_1].$ Let $z_0\in D$ such that it is in the domain of attraction of $\varphi.$ Then, for fix $\mu>0,$ the differential equation
\begin{equation}\label{eq1-noimpulse} \frac{dz}{dt} = \frac{F(z)}{\mu}  \end{equation}
with initial value $z(0,\mu)=z_0$ has a unique solution $z(t,\mu)$ since $F(z)\in C^1(D).$ Then rescale the time as  $t=\tau \mu$ and substitute $x=z-\varphi$ in \eqref{eq1-noimpulse} to get 
\begin{equation}\label{eq1-noimpulse2} 
	\frac{dx}{d\tau} = F(x+\varphi). 
\end{equation}
$x=0$ is an equilibrium of this differential equation. By condition (C1), equation \eqref{eq1-noimpulse2} has a positive definite function $V(x)$ whose derivative with respect to $\tau$ is negative definite and $V(0)=0.$  Hence, by the Lyapunuv second method, one can say that the zero solution of \eqref{eq1-noimpulse2} is uniformly asymptotically stable as $\tau \to \infty.$ Therefore, $\forall \, \varepsilon>0$ and for sufficiently small $\mu$ on $0<t\leq \theta_1$ one has $\|z(t,\mu)-\varphi\|<\varepsilon,$ i.e, 
\[\lim_{\mu \to 0}z(t,\mu)=\varphi \quad \text{for} \quad 0<t\leq \theta_1.\]
 Now, consider the next interval $(\theta_1,\theta_2].$ From  condition (C2), we have  \[\lim_{\mu \to 0}z(\theta_1+,\mu)= \lim_{\mu \to 0} \left\lbrace z(\theta_1,\mu)+\frac{I(z(\theta_1,\mu),\mu)}{\mu}\right\rbrace = \varphi.\]
It means that $z(\theta_1+,\mu)$ is in the domain of attraction of the root $\varphi.$ Repeating the same processes as for the previous interval, one obtains 
\[\lim_{\mu \to 0}z(t,\mu)=\varphi \quad \text{for} \quad \theta_1< t\leq \theta_2.\]
 Similarly, one can show that  $z(t,\mu)\to \varphi$ as $\mu \to 0$ for  $t\in (\theta_i,\theta_{i+1}],i=2,\dots,p-1$ and $t\in (\theta_p,T].$ As a result, limit \eqref{lim1} is true and the theorem is proved.

The convergence is not uniform at $t=0$ since $z(0,\mu)=z_0$ and $z_0\neq\varphi$ for all $\mu>0.$ We can say that the region of nonuniform convergence is  $O(\mu)$ thick, since for $t>0,$ $\|z(t,\mu)-\varphi\|$ can be made arbitrarily close to zero by choosing $\mu$ small enough. The interval of nonuniform convergence is called an initial layer. This theorem implies that there is a single initial layer.
 
\begin{figure}[H]
 	\centering
 	\includegraphics[scale=0.6]{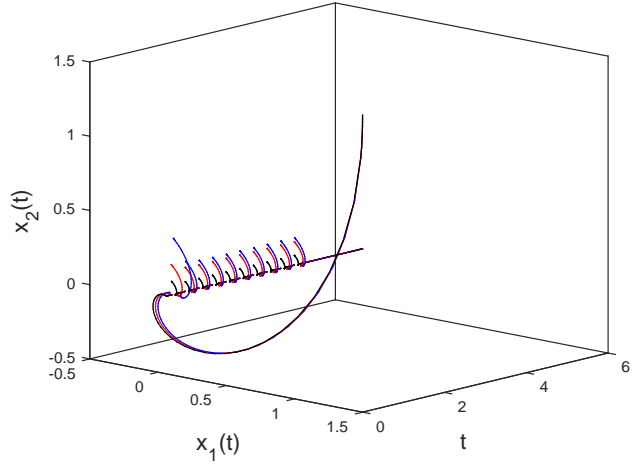}
 	\caption{Blue, red and black lines represents the solution of system \eqref{ex1} with initial value $(1.5,-1.5)$ for $\mu=0.07,\mu=0.05$ and $\mu=0.03,$ respectively.}
 	\label{ex1figure}
 \end{figure}

\textbf{Example. } Consider the system
\begin{equation}\label{ex1}
\begin{split}
&\mu \dot{x}_1=-x_1+x_2,\quad \mu \Delta x_1|_{t=\theta_i}=-2\mu x_1,\\
&\mu \dot{x}_2=-x_1-x_2,\quad \mu \Delta x_2|_{t=\theta_i}=\mu\sin(x_2^{1/3}+\mu),
\end{split}
\end{equation}
with initial value $(x_1(0,\mu),x_2(0,\mu)),$ where $\theta_i=i/3,i=1,2,\dots,10.$ Let us take $\mu=0$ in this system. Then 
\begin{equation*}
\begin{split}
&0=-x_1+x_2,\quad 0=0,\\
&0=-x_1-x_2,\quad 0=0.
\end{split}
\end{equation*}
and so $(x_1,x_2)=(0,0)$ is the root. Substitute $\tau=\frac{t}{\mu}$ into the differential equations part of \eqref{ex1} to obtain
\begin{equation}\label{ex1tau}
\begin{split}
&\frac{dx_1}{d\tau}=-x_1+x_2,\\
&\frac{dx_2}{d\tau}=-x_1-x_2,.
\end{split}
\end{equation}
We take  the  positive definite function $V(x_1,x_2)=x_1^2+x_2^2.$ Then
\[\frac{dV}{d\tau}=2x_1(-x_1+x_2)+2x_2(-x_1-x_2)=-2(x_1^2+x_2^2)=-2V.\]
 Hence, $V(x_1,x_2)$ has a negative definite derivative with respect to $\tau$ along \eqref{ex1tau}. Now, let us check the condition (C2). Denote $x=(x_1,x_2).$ Then 
 \[\lim_{(x,\mu)\to (0,0)}\frac{I(x,\mu)}{\mu}=0\]
 since $\lim_{(x_1,\mu)\to (0,0)}-2x_1=0$ and $\lim_{(x_2,\mu)\to (0,0)}\sin(x_2^{1/3}+\mu)=0.$
 Therefore, by Theorem \ref{thmlyapunov1}, if the initial value $(x_1(0,\mu),x_2(0,\mu))$ of \eqref{ex1} is in the domain of attraction of the root $(0,0),$ then solution $(x_1(t,\mu),x_2(t,\mu))$ of \eqref{ex1} tends to $(0,0)$ as $\mu \to 0$ for $0<t\leq T.$ It is clearly seen in  Figure \ref{ex1figure} that the solution of system \eqref{ex1} with initial $(1.5,-1.5)$ tends to $(0,0)$ as $\mu \to 0.$
% \begin{figure}[H]
% 	\centering
% 	\includegraphics[scale=0.6]{ex1}
% 	\caption{Blue and red lines represents the coordinates of  solution of system \eqref{ex1} with initial $(1.5,-1.5)$ for $\mu=0.3$ and $\mu=0.1,$ respectively.}
% 	\label{ex1figure}
% \end{figure}
%\textbf{Example 2.} Now, let us consider the following one dimensional nonlinear system 
%\begin{equation}\label{ex2}
%\begin{split}
%&\mu \dot{z}=-z-z^3,\\
%&\mu \Delta z|_{t=\theta_i}=z^2-z,
%\end{split}
%\end{equation}
%%%%%%%%%%%%%%%%%%%%%%%%%%%%%%
% Multi-Layers
%%%%%%%%%%%%%%%%%%%%%%%%%%%%%%
\subsection{Singularity with Multi-Layers}\label{smulti}
In the previous subsection, it is shown that there is a single initial layer. Using an impulse function, the convergence can be nonuniform near several points, that is to say that multi-layers emerge. These layers will occur on the neighborhoods of $t=0$ and $t=\theta_i,i=1,2,\dots,p.$

Again, we consider  system \eqref{eq1} with the same properties. In addition, we need the following condition
 \begin{itemize}
 \item[(C3)]
 \[\lim_{(z,\mu)\to (\varphi,0)}\frac{I(z,\mu)}{\mu}=I_0\neq 0\] and assume that $\varphi+I_0$ is in the domain of attraction of the root $\varphi.$
 \end{itemize}
 By the virtue of this condition, after the each impulse moment, the difference $\|z(\theta_i+,\mu)-\varphi\|$ does not go to zero as $\mu \to 0.$ Hence, convergence is not uniform.
 
\begin{theorem}\label{multilayer}
Suppose that conditions (C1) and (C3) hold. If the initial value $z_0$ is located in the domain of attraction of the root $\varphi,$ then  solution $z(t,\mu)$ of \eqref{eq1} with $z(0,\mu)=z_0$ exists on $[0,T]$ and the limit
\begin{equation}\label{lim2}
	\lim_{\mu \to 0}z(t,\mu)=\varphi
\end{equation}
is true for $ t\in  \bigcup\limits_{i=0}^{p-1}(\theta_i,\theta_{i+1}]\cup(\theta_p,T],$
where $\theta_0=0.$
\end{theorem}
\textbf{Proof.}
Proof is similar to the proof of Theorem \ref{thmlyapunov1} with the exception that singularity with multi-layers appears near $t=0$ and $t=\theta_i,i=1,2,\dots,p.$

By condition (C3), after the each discontinuity moment $\theta_i,$ the solution $z(t,\mu)$ is not close to the root $\varphi.$ In other words, the difference  $\|z(\theta_i+,\mu)-\varphi\|$ cannot be arbitrarily small as $\mu \to 0.$ Hence, one can obtain  multi-layers up to the number $p+1.$

 Let us illustrate the theorem with the following example.
 
\begin{figure}[H]
\centering
\includegraphics[scale=0.5]{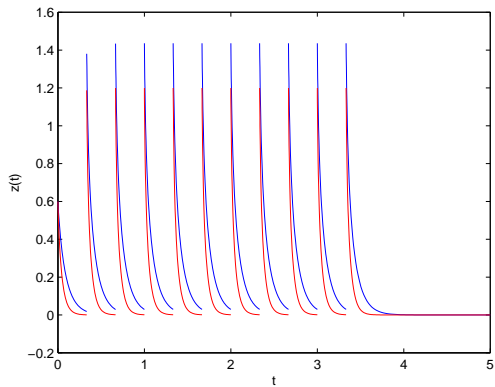}
\caption{Solution $z(t,\mu)$ of system \eqref{ex2} with initial value $z(0,\mu)=0.6$ for different values of parameter $\mu.$ Blue and red line represent for $\mu=0.1,\mu=0.05,$ respectively. It is seen that at $t=0$ and at each $\theta_i,i=1,2,...,10,$ the convergence is nonuniform, i.e., multi-layers exist.}
\label{fig:multilayer}
\end{figure}

\textbf{Example.} Consider the following impulsive differential equation with small parameter:
\begin{equation}\label{ex2}
\begin{split}
	&\mu \dot{z}=-z-z^3, \\
	&\mu \Delta z|_{t=\theta_i}=\mu z^{1/3}+\sin\mu+0.1\mu,
\end{split}
\end{equation}
where $\theta_i=i/3,i=1,2,...,10.$
Let us take $\mu=0$ in this system. Then we have the algebraic equation $0=-z-z^3.$ It has solution $z=0.$ Now, introduce $t=\tau \mu$ in the first equation of \eqref{ex2} to  obtain 
\begin{equation}\label{ex2tau}
\frac{dz}{d\tau}=-z-z^3
\end{equation}
Using the Lyapunov function $V(z)=z^2,$ it can be shown that $z=0$ is a uniformly asymptotic stable solution of \eqref{ex2tau}. Moreover, condition (C3) is satisfied since
\[\lim_{(z,\mu)\to (0,0)}\frac{\mu z^{1/3}+\sin\mu+\mu0.1}{\mu}=1.1.\]
Choose the initial value $z(0,\mu)=0.6.$ Then the solution $z(t,\mu)$ of system \eqref{ex2} with this initial value has multi-layers at $t=0$ and  $t=\theta_i+,i=1,2,\dots,10.$ Clearly, in Figure \ref{fig:multilayer}, it can be seen that multi-layers occur.
 Let us generalize the Theorem \ref{multilayer}. Consider the following impulsive system
  \begin{equation} \label{multigenel}
  \begin{split}
  &\mu \frac{dz}{dt} = F(z) , \\
  &\mu \Delta z|_{t=\theta_{i}}=I(z,\mu), \\
  &\mu \Delta z|_{t=\tau_{j}^i}=J(z,\mu),
  \end{split}
  \end{equation}
 where  the the impulse moments $\tau_j^i,j=1,2,\dots,p_j$ are such that $\theta_i<\tau_1^i<\tau_2^i<\dots<\tau_{p_j}^i, i=1,2,\dots,p-1$ and $\theta_p<\tau_1^p<\tau_2^p<\dots<\tau_{p_j}^p<T.$ Assume the following condition holds for \eqref{multigenel}
 \begin{itemize}
 	\item[(C4)] \[ \lim_{(z,\mu)\to (\varphi,0)}\frac{J(z,\mu)}{\mu}= 0,\]
 \end{itemize}
 Now, we can assert the following theorem.
 \begin{theorem}\label{multilayer2}
 	Suppose that conditions (C1), (C3) and (C4) hold. If the initial value $z_0$ is located in the domain of attraction of the root $\varphi,$ then  solution $z(t,\mu)$ of \eqref{multigenel} with $z(0,\mu)=z_0$ exists on $[0,T]$ and the limit
 	\begin{equation}
 	\lim_{\mu \to 0}z(t,\mu)=\varphi
 	\end{equation}
 	is true for $ t\in  \bigcup\limits_{i=0}^{p-1}(\theta_i,\theta_{i+1}]\cup(\theta_p,T],$
 	where $\theta_0=0.$
 \end{theorem}
%%%%%%%%%%%%%%%%%%%%%%%%%%%%%%%%%%%%%%%%%%%%%%%%%%%%%%%%
%%%
%%%%%%%%%%%%%%%%%%%%%%%%%%%%%%%%%%%%%%%%%%%%%%%%%%%%%%%%
\section{Main Result}
Now, we turn to  main problem \eqref{eq1-general}.
\subsection{Singularity with  a Single Layer}
Define the initial conditions (for simplicity, we set $t_0 = 0$.)
\begin{equation} \label{eq1-initials}
z(0,\mu)=z^{0}, y(0,\mu)=0,
\end{equation}
where $z^0$ and $y^0$ will be assumed to be independent of $\mu$, and let us investigate the solution $z (t, \mu)$, $y (t, \mu)$ of \eqref{eq1-general} and \eqref{eq1-initials} on segment $0\leq t\leq T$.

In system \eqref{eq1-general}, take $\mu=0$, then we obtain
\begin{equation} \label{eq1-degenerate}
\begin{split}
&0 = F(\bar{z},\bar{y},t), \quad 0=I(\bar{z},\bar{y},0),\\
&\frac{d\bar{y}}{dt} = f(\bar{z},\bar{y},t), \quad \Delta \bar{y}|_{t=\eta_{j}}=J(\bar{z},\bar{y})
\end{split}
\end{equation} 
which we call as degenerate system  due to the fact that its order is less than the order of \eqref{eq1-general}. Therefore, for the system \eqref{eq1-degenerate} the number of initial conditions must be set less than the number of initial conditions for \eqref{eq1-general}. We naturally insert the initial condition for $y$, i.e., put
\begin{equation} \label{initial-y0}
\bar{y}(0)=y^0,
\end{equation} 
and drop the initial condition for $z$.
Now, the question is that whether there will be a solution $z(t,\mu)$ and $y(t,\mu)$ of problem \eqref{eq1-general}, \eqref{eq1-initials} for small $\mu$  which is close to the solution $\bar{z}(t), \bar{y}(t)$ of the degenerate problem \eqref{eq1-degenerate}, \eqref{initial-y0}.

To solve system \eqref{eq1-degenerate}, it is necessary to find $\bar{z}$ from $0=F(\bar{z},\bar{y},t)$ and $0=I(\bar{z},\bar{y},0).$ Then choose one of the root $\bar{z}=\varphi(\bar{y},t)$ such that $0=F(\varphi(\bar{y},t),\bar{y},t)$ and $0=I(\varphi(\bar{y},\theta_i),\bar{y},0),$ and substitute into \eqref{eq1-degenerate} with initial value \eqref{initial-y0} to obtain
\begin{equation}	\label{ybar}
\begin{split}
&\frac{d\bar{y}}{dt} = f(\varphi(\bar{y},t),\bar{y},t), \quad \Delta \bar{y}|_{t=\eta_{j}}=J(\varphi(\bar{y},t),\bar{y}), \\ 
&\bar{y}(0)=y^0.
\end{split}
\end{equation}
We need  the following conditions in this section:
\begin{itemize}
	\item[A1.] The functions $F(z,y,t), f(y,z,t),$ and $J(z,y)$ are continuous in some domain $H=\{(y,t)\in\bar{N}=\{0\leq t\leq T, \|y\|\leq c\}, \|z\|<d \}$, $I(z,y,\mu)$ is continuous  in $H\times [0,1]$ and they are Lipschitz  continuous with respect to $z$ and $y$.
	\item[A2.] Algebraic equations $0=F(z,y,t)$ and $0=I(z,y,0)$ have a root $z=\varphi(y,t)$ such that $F(\varphi(y,t),y,t)=0$  and  $I(\varphi(y(\theta_i),\theta_i),y(\theta_i),0)=0,i=1,2,\dots,p,$ in  domain $\bar{N}$  such that:
	\begin{enumerate}
		\item  $\varphi(y,t)$ is a continuous function in $\bar{N}$,
		\item ($\varphi(y,t),y,t)$ $\in$ H,
		\item The root $\varphi(y,t)$ is isolated in $\bar{N}$, i.e., $\exists \, \epsilon>0$: $F(z,y,t)\neq 0$ and/or $I(z,y,\mu)\neq 0$ for $0<\|z-\varphi(y,t)\|<\epsilon$, $(y,t) \in \bar{N}.$
	\end{enumerate}
	\item[A3.]
	\begin{enumerate}
		\item System \eqref{ybar} has a unique solution $\bar{y}(t)$ on $0\leq t \leq T$, and $(\bar{y}(t),t) \in \bar{N}$ for $0\leq t \leq T$. Moreover, $f(\varphi(y,t),y,t)$ and $J(\varphi(y,t),y)$ are	Lipschitz with respect to $y \in \bar{N}$.
		\item $\varphi(\bar{y}(\eta_j+),\eta_j+)=\varphi(\bar{y}(\eta_j),\eta_j),j=1,2,\dots,k.$
	\end{enumerate} 
\end{itemize}
Now, setting $x=z-\varphi$ and $t=\tau \mu,$ we introduce the  system
\begin{equation} \label{eq1-adjoint}
\frac{dx}{d\tau} = F(x+\varphi(y,t),y,t), \quad \tau \geq 0,
\end{equation}
where  $y$ and $t$ are considered as parameters, $x=0$ is an isolated stationary point of \eqref{eq1-adjoint} for $(y,t) \in \bar{D}$.	
\begin{itemize}		
	\item[A4.]Suppose that  there is a positive definite function $V(x,y,\tau)$ whose derivative with respect to $\tau$ along the system \eqref{eq1-adjoint} is negative definite in the region $H.$
\end{itemize}
Consider adjoint system 
\begin{equation} \label{eq1-adjoint-thm}
\frac{d\tilde{z}}{d\tau} = F(\tilde{z},y^0,0),% 0=I(\tilde{z},y^0,0) 
\quad\tau \geq 0, \\
\end{equation}
with initial condition 
\begin{equation}
\tilde{z}(0)=z^0.
\label{z0-thm}
\end{equation}

Since $z^0$ maybe, in general, far from stationary point $\varphi(y^0,0)$, then the solution $\tilde{z}(\tau)$ of equations \eqref{eq1-adjoint-thm} and \eqref{z0-thm} need not tend to $\varphi(y^0,0)$ as $\tau \rightarrow \infty$. Assume also that
\begin{itemize}
	\item[A5.] the solution $\tilde{z}(\tau)$ of equations \eqref{eq1-adjoint-thm} and \eqref{z0-thm} satisfies the conditions
	\begin{enumerate}
		\item $\tilde{z}(\tau) \to \varphi(y^0,0)$ as $\tau \to \infty$,
		\item $(\tilde{z}(\tau),y^0,0) \in H$ for $\tau \geq 0.$
	\end{enumerate}
\end{itemize}
In this case, $z^0$ is said to belong to the basin of attraction of the stationary point $\tilde{z}=\varphi(y^0,0)$. By virtue of the asymptotic stability of this point all points near it will belong to its basin of attraction.
\begin{itemize}
	\item[A6.] Assume also 
	\[\lim_{(z,y,\mu)\to (\varphi,y^0,0)}\frac{I(z,y,\mu)}{\mu}=0.\]
\end{itemize}
Now, we state and  prove the modified Tikhonov Theorem.
\begin{theorem} \label{thm-general}
	Suppose that  conditions $A1-A6$ hold. Then, for sufficiently small $\mu,$ solutions $z(t,\mu)$ and $y(t,\mu)$ of problem \eqref{eq1-general} with initial conditions \eqref{eq1-initials} exist on $0\leq t \leq T$, are unique, and satisfy
	\begin{equation}\label{thm-general-result1}
	\lim_{\mu \to 0}y(t,\mu)=\bar{y}(t) \quad \textit{for}\quad 0 \leq t \leq T
	\end{equation}
	and
	\begin{equation}\label{thm-general-result2}
	\lim_{\mu \to 0}z(t,\mu)=\bar{z}(t)=\varphi(\bar{y}(t),t) \quad \textit{for} \quad 0 < t \leq T.
	\end{equation}
\end{theorem}
Before proving this theorem, we will consider the following  auxiliary system:
\begin{equation} \label{eq1-general2}
\begin{split}
\mu \frac{dz}{dt} &= F(z,y,t) ,\\
\frac{dy}{dt} &= f(z,y,t) , \quad \Delta y|_{t=\eta_{j}}=J(z,y),
\end{split}
\end{equation} 
where this system has same properties as \eqref{eq1-general}.

In system \eqref{eq1-general2}, take $\mu=0$, then we obtain
\begin{equation} \label{eq1-degenerate2}
\begin{split}
&0 = F(\bar{z},\bar{y},t),\\
&\frac{d\bar{y}}{dt} = f(\bar{z},\bar{y},t), \quad \Delta \bar{y}|_{t=\eta_{j}}=J(\bar{z},\bar{y})
\end{split}
\end{equation} 
which is degenerate system  of \eqref{eq1-general2}.

To solve system \eqref{eq1-degenerate2}, it is necessary to find $\bar{z}$ from $0=F(\bar{z},\bar{y},t).$ Then choose one of the root $\bar{z}=\varphi(\bar{y},t)$ and substitute into \eqref{eq1-degenerate2} with initial value \eqref{initial-y0} to obtain
\begin{equation}	\label{ybar2}
\begin{split}
&\frac{d\bar{y}}{dt} = f(\varphi(\bar{y},t),\bar{y},t), \quad \Delta \bar{y}|_{t=\eta_{j}}=J(\varphi(\bar{y},t),\bar{y}), \\ 
&\bar{y}(0)=y^0.
\end{split}
\end{equation}

Now, introduce the adjoint system
\begin{equation} \label{eq1-adjoint2}
\frac{d\tilde{z}}{d\tau} = F(\tilde{z},y,t) \quad \tau \geq 0, \\ 
\end{equation}
where $y$ and $t$ are considered as parameters, $\tilde{z}=\varphi(y,t)$ is an isolated stationary point of \eqref{eq1-adjoint2} for $(y,t) \in \bar{N}$. 

Suppose that
\begin{itemize}		
	\item[B.] the stationary point $\tilde{z}=\varphi(y,t)$ of \eqref{eq1-adjoint} is  uniformly asymptotically stable with respect to $(y,t) \in \bar{N}$, i.e. $\forall \varepsilon>0$ $\exists\, \delta(\varepsilon)>0$ such that if $\|\tilde{z}(0)-\varphi(y,t)\|<\delta(\varepsilon)$ then
	$ \|\tilde{z}(\tau)-\varphi(y,t)\|<\varepsilon$ and $ \tilde{z}(\tau)\to \varphi(y,t)$ as $ \tau \to\infty .$
\end{itemize}	
If this condition is true, then the root $\tilde{z}=\varphi(y,t)$ is said to be stable in $\bar{N}$.

%Consider adjoint system \eqref{eq1-adjoint} at $t=0, y=y^0$
%\begin{equation} \label{eq1-adjoint0}
%\frac{d\tilde{z}}{d\tau} = F(\tilde{z},y^0,0) \quad \tau \geq 0, \\
%\end{equation}
%with initial condition \eqref{z0-thm}.
%Since $z^0$ maybe, in general, far from stationary point $\varphi(y^0,0)$, then the solution $\tilde{z}(\tau)$ of equations \eqref{eq1-adjoint0} and \eqref{z0-thm} need not tend to $\varphi(y^0,0)$ as $\tau \rightarrow \infty$. Assume also that
%\begin{itemize}
%	\item[B2.] the solution $\tilde{z}(\tau)$ of equations \eqref{eq1-adjoint0} and \eqref{z0-thm} satisfies the conditions
%	\begin{enumerate}
%		\item $\tilde{z}(\tau) \to \varphi(y^0,0)$ as $\tau \to \infty$,
%		\item $(\tilde{z}(\tau),y^0,0) \in H$ for $\tau \geq 0.$
%	\end{enumerate}
%\end{itemize}
%Thus, $z^0$  belongs to the basin of attraction of the stationary point $\tilde{z}=\varphi(y^0,0)$. By virtue of the asymptotic stability of this point all points near it will belong to its basin of attraction.
\begin{lemma} \label{lemma}
	Suppose that for system \eqref{eq1-general2}  conditions A1-A3,A5 and B are true, then, for sufficiently small $\mu,$ solutions $z(t,\mu)$ and $y(t,\mu)$ of problem \eqref{eq1-general2} with initial conditions \eqref{eq1-initials} exist on $0\leq t \leq T$, are unique, and satisfy
	\begin{equation}\label{lemma-result1}
	\lim_{\mu \to 0}y(t,\mu)=\bar{y}(t) \quad \textit{for}\quad 0 \leq t \leq T
	\end{equation}
	and
	\begin{equation}\label{lemma-result2}
	\lim_{\mu \to 0}z(t,\mu)=\bar{z}(t)=\varphi(\bar{y}(t),t) \quad \textit{for} \quad 0 < t \leq T.
	\end{equation}
\end{lemma}
\textbf{Proof. } First, consider the interval $[0,\eta_1].$ On this interval, Lemma \ref{lemma} is type of  Tikhonov Theorem \cite[Theorem 2.1]{vasileva} and all conditions are satisfied. Therefore, by \cite[Theorem 2.1]{vasileva}, for sufficiently small $\mu,$ solutions $z(t,\mu),y(t,\mu)$ of \eqref{eq1-general} and \eqref{eq1-initials} exist on $[0,\eta_1]$ and satisfies 
\begin{equation}
\begin{split}\label{prooflimit}
&\lim_{\mu \to 0}y(t,\mu)=\bar{y}(t) \quad \textit{for}\quad 0 \leq t \leq \eta_1,\\ &\lim_{\mu \to 0}z(t,\mu)=\bar{z}(t)=\varphi(\bar{y}(t),t) \quad \textit{for} \quad 0 < t \leq \eta_1.
\end{split}
\end{equation}
Now, consider the second interval $(\eta_1,\eta_2].$ For this interval the initial values are $z_1=z(\eta_1+,\mu),y_1=y(\eta_1+,\mu).$ Since $\lim_{\mu \to 0}y(\eta_1,\mu)=\bar{y}(\eta_1)$ and $\lim_{\mu \to 0}z(\eta_1,\mu)=\varphi(\bar{y}(\eta_1),\eta_1),$ $z_1$ is in the the basin of attraction of $\varphi(\bar{y}(t),t)$ and $y_1\in N.$ Again, all conditions of Tikhonov Theorem are satisfied and by \cite[Theorem 2.1]{vasileva}
\begin{equation*}
\begin{split}
&\lim_{\mu \to 0}y(t,\mu)=\bar{y}(t) \quad \textit{for}\quad \eta_1 < t \leq \eta_2,\\ &\lim_{\mu \to 0}z(t,\mu)=\bar{z}(t)=\varphi(\bar{y}(t),t) \quad \textit{for} \quad \eta_1 < t \leq \eta_2.
\end{split}
\end{equation*}
Similarly, for the next intervals $(\eta_i,\eta_{i+1}],i=2,3,\dots,k-1,$ and $(\eta_k,T]$ one can show that as $\mu \to 0,$ $\lim_{\mu \to 0}y(t,\mu)=\bar{y}(t)$ and $\lim_{\mu \to 0}z(t,\mu)=\bar{z}(t)=\varphi(\bar{y}(t),t).$ Lemma is proved.

\textbf{Remark. } At discontinuity moments $\eta_j,j=1,2,\dots,k,$ layers do not emerge. This is because, $\varphi(\bar{y}(t),t)$ is a continuous function and $\lim_{\mu \to 0}z(\eta_j+,\mu)=\varphi(\bar{y}(\eta_j),\eta_j)=\bar{z}(\eta_j),j=1,2,\dots,k.$
%since $z(t,\mu)$ does not have any impulse function once it is in $\varepsilon-$neighborhood of $\varphi(y,t),$ it never leaves this neighborhood. Hence, for $0< t \leq T$ whenever $(y(t,\mu),t)\in N$  we can write 
%\begin{equation}\label{gamma}
%z(t,\mu)=\varphi(y,t)+\gamma(t,\mu),
%\end{equation}
%where $\gamma(t,\mu)$ is a continuous function satisfying $\|\gamma(t,\mu)\|<\varepsilon.$ Substitute \eqref{gamma} into the second line in \eqref{eq1-general2} to get
%\begin{equation}	\label{ygamma}
%\begin{split}
%&\frac{dy}{dt} = f(\varphi(y,t)+\gamma(t,\mu),y,t), \quad \Delta y|_{t=\eta_{j}}=J(\varphi(y,t)+\gamma(t,\mu),y), \\ 
%&y|_{t=t_0}=y^0+\omega(\mu),
%\end{split}
%\end{equation}
%where $t_0=t_0(\mu)<\varepsilon$ and $\omega(\mu) \to 0$ as $\mu \to 0.$ System \eqref{ygamma} is a regularly perturbed system with respect to \eqref{ybar2}. Thus, from the regular perturbation theorem for discontinuous dynamical systems \cite{Akhmet2010}, for sufficiently small $\mu,$ one has 
%\[\|y(t,\mu)-\bar{y}(t)\|<\varepsilon \text{ for } t_0\leq t \leq T.\]
%From this inequality and the first limit in \eqref{prooflimit}, we obtain \eqref{lemma-result1}. Therefore, lemma is proved.

\textbf{Proof of Theorem \ref{thm-general}. }  Consider the interval $[0,\theta_1].$ Hence, on this interval, Theorem \ref{thm-general} is type of  Lemma \ref{lemma}. Condition A4.  is corresponding to the assumption that uniform asymptomatic stability of the root $\varphi(y,\tau\mu)$ as $\tau \to \infty,$ i.e. condition B is satisfied. Obviously, all conditions of the lemma are true. Consequently, for sufficiently small $\mu,$ solutions $z(t,\mu),y(t,\mu)$ of \eqref{eq1-general2} and \eqref{eq1-initials} exist and satisfies 
\begin{equation}
\begin{split}\label{prooflimit2}
&\lim_{\mu \to 0}y(t,\mu)=\bar{y}(t) \quad \textit{for}\quad 0 \leq t \leq \theta_1,\\ &\lim_{\mu \to 0}z(t,\mu)=\bar{z}(t)=\varphi(\bar{y}(t),t) \quad \textit{for} \quad 0 < t \leq \theta_1.
\end{split}
\end{equation}
Now, consider the next interval $(\theta_1,\theta_2].$ Condition A6 implies that
\[ \lim_{\mu\to 0}z(\theta_1+,\mu)=	\lim_{\mu \to 0}\left\lbrace z(\theta_1,\mu)+ \frac{I(z(\theta_1,\mu),y(\theta_1,\mu),\mu)}{\mu}\right\rbrace  =\varphi. \]  Hence, condition A5 is true.  Repeating the same processes as for the previous interval, one can demonstrates that $z(t,\mu)\to \varphi(\bar{y}(t),t)$ and $y(t,\mu)\to \bar{y}(t)$ as $\mu \to 0$ for $(\theta_1,\theta_2].$ Thus, recurrently it can be proven that for $t\in (\theta_i,\theta_{i+1}],i=1,2,\dots,p-1$ and $t\in (\theta_p,T]$ it is true that $z(t,\mu)\to \varphi(\bar{y}(t),t)$ and $y(t,\mu)\to \bar{y}(t)$ as $\mu \to 0.$ Therefore limits \eqref{thm-general-result1} and \eqref{thm-general-result2} are true. Theorem is proved.

\begin{figure}[H]
	\centering
	\includegraphics[scale=0.75]{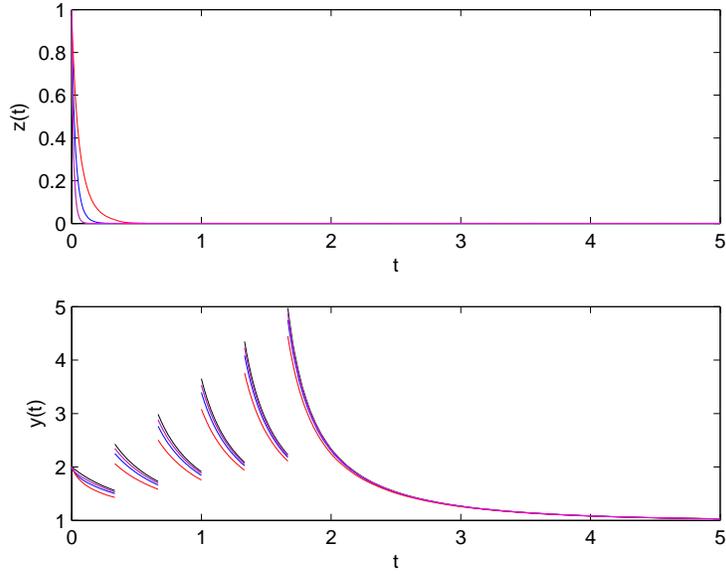}
	\caption{Black,magenta, blue and red lines are the coordinates of system \eqref{ex11} with initial values $z(0,\mu)=1$ and $y(0,\mu)=2$ for different values of $\mu:0,0.05,0.1,0.2,$ respectively.}
	\label{lemmaex}
\end{figure}

\textbf{Example for Lemma \ref{lemma}. } Consider the system
\begin{equation}\label{ex11}
\begin{split}
&\mu \frac{dz}{dt}=z(1-z-2y),\\
& \frac{dy}{dt}=y(1-2z-y), \quad \Delta y|_{t=\eta_j}=y^2-y+z,
\end{split}
\end{equation}
with initial conditions $z(0,\mu)=1$ and $y(0,\mu)=2,$ where $\eta_j=j/3,j=1,2,\dots,5.$ Let us take $\mu=0$ in this problem. Then, the first equation becomes $0=z(1-z-2y).$ It has the solutions $z=0$ and $z=1-2y.$ Consider the zero solution $z=0.$ Now, we check the conditions of Lemma \ref{lemma}. \[\frac{\partial}{\partial z}z(1-z-2y)|_{z=0}={1-2y}<0\] if $y>1/2.$ Therefore, if $y>1/2,$ $z=0$ is uniformly asymptotically stable. Substitute $z=0$ into the second line of \eqref{ex11} to obtain 
\begin{equation}\label{ex1-degenerate}
\frac{d\bar{y}}{dt}=\bar{y}(1-\bar{y}), \quad \Delta \bar{y}|_{t=\eta_j}=\bar{y}^2-\bar{y},
\end{equation}
with initial value $\bar{y}(0)=2.$ This system has a unique solution $\bar{y}(t)$. Thus, by Lemma \ref{lemma}, solutions $z(t,\mu),y(t,\mu)$ of \eqref{ex11} with $z(0,\mu)=1$ and $y(0,\mu)=2$ tends to $0,\bar{y}(t),$ respectively, as $\mu \to 0$ for $0<t\leq T.$ Obviously, in Figure \ref{lemmaex}, it can be seen that when $\mu$ decreases to zero, solutions $z(t,\mu),y(t,\mu)$ approaches to  $0,\bar{y}(t),$ respectively.
%%%%%%%%%%%%%%%%%%%%%%%%%%%%%%
% Multi-Layers
%%%%%%%%%%%%%%%%%%%%%%%%%%%%%%
\subsection{Singularity with Multi-Layers}
In the previous subsection, we have shown that the convergence is not uniform at $t=0.$ That is,  an initial layer is obtained by Tikhonov Theorem. To get  multi-layers by Tikhonov Theorem we need another condition for the impulse function. These layers will occur on the neighborhoods of $t=0$ and $t=\theta_i,i=1,2,\dots,p.$

Again, we consider  system \eqref{eq1-general} with the same properties. In addition, we need the following condition
\begin{itemize}
	\item[A7.]
	\[\lim_{(z,y,\mu)\to (\varphi,\bar{y},0)}\frac{I(z,y,\mu)}{\mu}=I_0\neq 0\] and assume that $\varphi(\bar{y}(\theta_i),\theta_i)+I_0,i=1,2,\dots,p,$ is in the basin of attraction of $\varphi(\bar{y}(t),t).$
\end{itemize}
This condition implies that after each impulse moment, the difference $\|z(\theta_i+,\mu)-\varphi\|$ does not go to zero as $\mu \to 0.$ Hence, convergence is not uniform.
\begin{theorem} 
	Suppose that  conditions A1-A5 and A7 hold. Then, for sufficiently small $\mu,$ solutions $z(t,\mu)$ and $y(t,\mu)$ of problem \eqref{eq1-general} with initial conditions \eqref{eq1-initials} exist on $0\leq t \leq T$, are unique, and satisfy
	\[
	\lim_{\mu \to 0}y(t,\mu)=\bar{y}(t) \quad \textit{for}\quad 0 \leq t \leq T
	\]
	and
	\[
	\lim_{\mu \to 0}z(t,\mu)=\bar{z}(t)=\varphi(\bar{y}(t),t)
	\]
	is true for $ t\in  \bigcup\limits_{i=0}^{p-1}(\theta_i,\theta_{i+1}]\cup(\theta_p,T],$
	where $\theta_0=0.$
\end{theorem}
\textbf{Proof.}
Proof is similar to the proof of Theorem \ref{thm-general} with the exception that singularity with multi-layers appears near $t=0$ and $t=\theta_i,i=1,2,\dots,p.$

Now, let us generalize this theorem. Consider the following impulsive system
\begin{equation} \label{multilayer3}
\begin{split}
\mu \frac{dz}{dt} &= F(z,y,t) ,\quad \mu \Delta z|_{t=\theta_{i}}=I(z,y,\mu)\quad \mu \Delta z|_{t=\tau_{j}^i}=\tilde{J}(z,y,\mu)\\
\frac{dy}{dt} &= f(z,y,t) , \quad \Delta y|_{t=\eta_{j}}=J(z,y),
\end{split}
\end{equation} 
where $\tau_{j}^i$ is defined in Subsection \ref{smulti}. Additionally, we need the condition 
\begin{itemize}
	\item[A8.]
	\[\lim_{(z,y,\mu)\to (\varphi,\bar{y},0)}\frac{J(z,y,\mu)}{\mu}=0.\]
\end{itemize}
Now we can assert  our theorem.
\begin{theorem} 
	Suppose that  conditions A1-A5 and A7-A8 hold. Then, for sufficiently small $\mu,$ solutions $z(t,\mu)$ and $y(t,\mu)$ of problem \eqref{multilayer3} with initial conditions \eqref{eq1-initials} exist on $0\leq t \leq T$, are unique, and satisfy
	\[
	\lim_{\mu \to 0}y(t,\mu)=\bar{y}(t) \quad \textit{for}\quad 0 \leq t \leq T
	\]
	and
	\[
	\lim_{\mu \to 0}z(t,\mu)=\bar{z}(t)=\varphi(\bar{y}(t),t)
	\]
	is true for $ t\in  \bigcup\limits_{i=0}^{p-1}(\theta_i,\theta_{i+1}]\cup(\theta_p,T],$
	where $\theta_0=0.$
\end{theorem}
\section{Conclusion} In this manuscript, we have introduced a new type of singular impulsive differential equation model. In this model,   Lyapunov second method is used to show the stability in the rescaled time. Then some illustrative examples with simulations are given to support the theoretical results.

The main novelty of this research is that singularity in the impulsive part of the systems can be treated through perturbation methods.

In the book of Baionov and Covachev \cite{bainov}, and several papers cited in the book, they considered singular impulsive systems with small parameter involved only in the differential equations of the systems, but not in the impulsive equations of them. This is why, we insert a small parameter into the the impulse  equation such that the singularity concept has been significantly extended for discontinuous dynamics.
\section*{References}

\bibliography{mybibfile}

\end{document}